\documentclass{ws-procs9x6}

\newcommand{\BEQ}{\begin{equation}}     % Gleichungen Anfang ..
\newcommand{\BEA}{\begin{eqnarray}}
\newcommand{\EEQ}{\end{equation}}       % .. und Ende
\newcommand{\EEA}{\end{eqnarray}}
\newcommand{\eps}{\varepsilon}          % epsilon
\newcommand{\g}{{\mathfrak{g}}}
\newcommand{\so}{{\mathfrak{so}}}
\newcommand{\h}{{\mathfrak{h}}}

\newcommand{\p}{{\mathfrak{p}}}
\newcommand{\sv}{{\mathfrak{sv}}}

\newcommand{\slin}{{\mathfrak{sl}}}

\newcommand{\alt}{{\mathfrak{alt}}}
\newcommand{\sch}{{\mathfrak{sch}}}
\newcommand{\conf}{{\mathfrak{conf}}}

\newcommand{\eop}{$\Box$}		% fin logique des demonstrations

\newcommand{\al}{\alpha}

\newcommand{\half}{\frac{1}{2}}

\newcommand{\R}{\mathbb{R}}		% ensembles reel, complexe, entier

\newcommand{\Z}{\mathbb{Z}}
\newcommand{\D}{{\rm d}}                % gerades d fuer Ableitungen
\newcommand{\II}{{\rm i}}               % gerades i fuer komplexe Einheit
\newcommand{\rar}{\rightarrow}          % Pfeil nach rechts
\begin{document}
\title{On the dynamical symmetric algebra of ageing: Lie structure, representations and Appell systems}
\author{Malte HENKEL
\thanks{LPM (CNRS UMR 7556), Universit\'e Henri Poincar\'e, BP 239, 
F - 54506 Vandoeuvre-l\`es-Nancy, France}, 
Ren\'e SCHOTT\thanks{IECN and LORIA, Universit\'e Henri Poincar\'e, BP 239, 
F - 54506 Vandoeuvre-l\`es-Nancy, France }, 
Stoimen STOIMENOV
\thanks{LPM (CNRS UMR 7556), Universit\'e Henri Poincar\'e, BP 239, 
F - 54506 Vandoeuvre-l\`es-Nancy, France and Institute of Nuclear Research and Nuclear Energy, Bulgarian Academy of Sciences, 1784 Sofia, Bulgaria},\\ 
J\'er\'emie UNTERBERGER
\thanks{IECN, Universit\'e Henri Poincar\'e, BP 239, \
F - 54506 Vandoeuvre-l\`es-Nancy, 
France}}
\date{}
\maketitle
\abstracts{The study of ageing phenomena leads to the investigation of a maximal parabolic subalgebra of $\conf_3$ which we call $\alt$. 
We investigate its Lie structure, prove some results concerning its representations and characterize the related Appell systems.}
\section{Introduction}
\begin{figure}[t]
\centerline{\epsfxsize=1.5in\ \epsfbox{
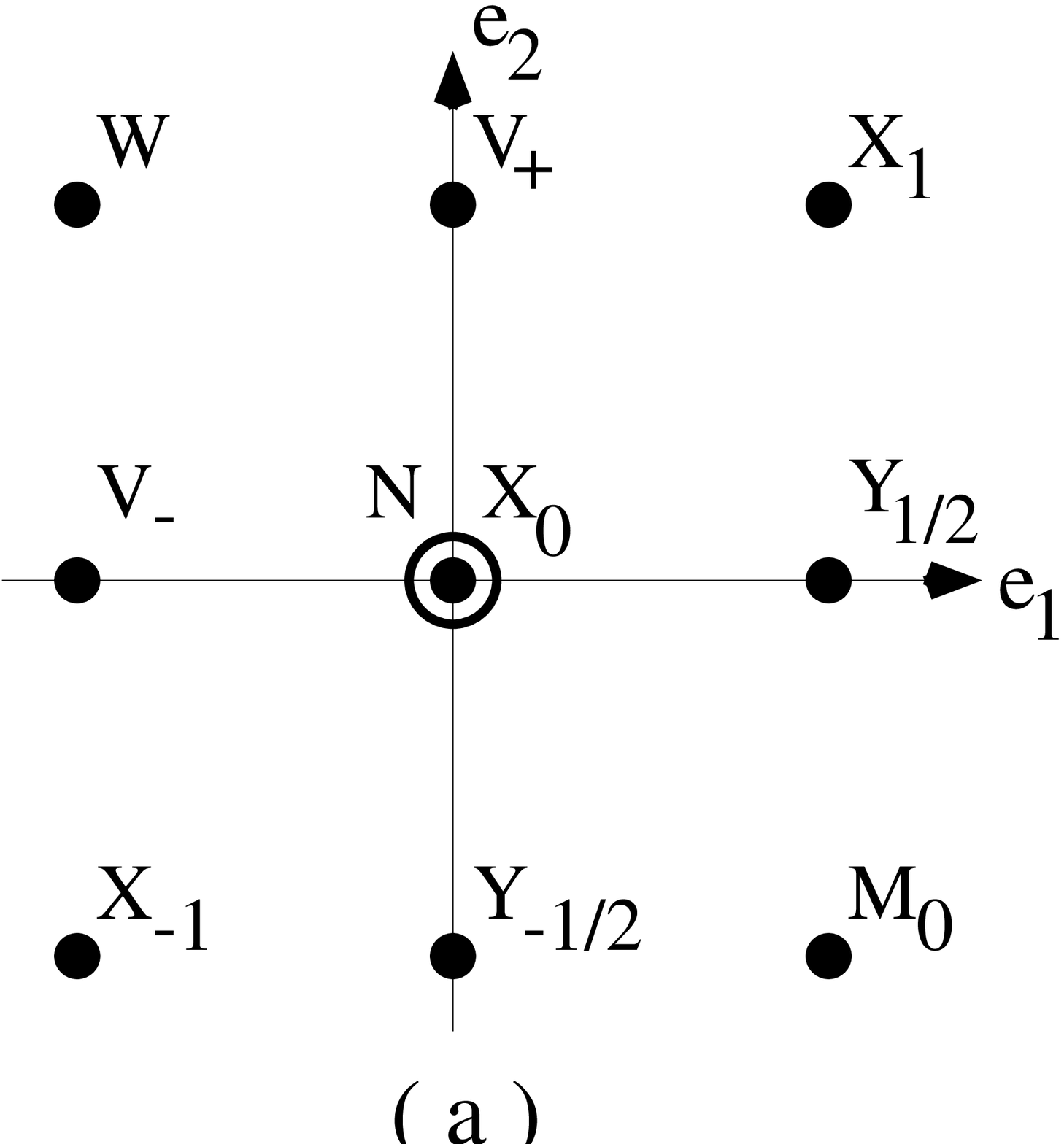} ~
\epsfxsize=1.5in\epsfbox{
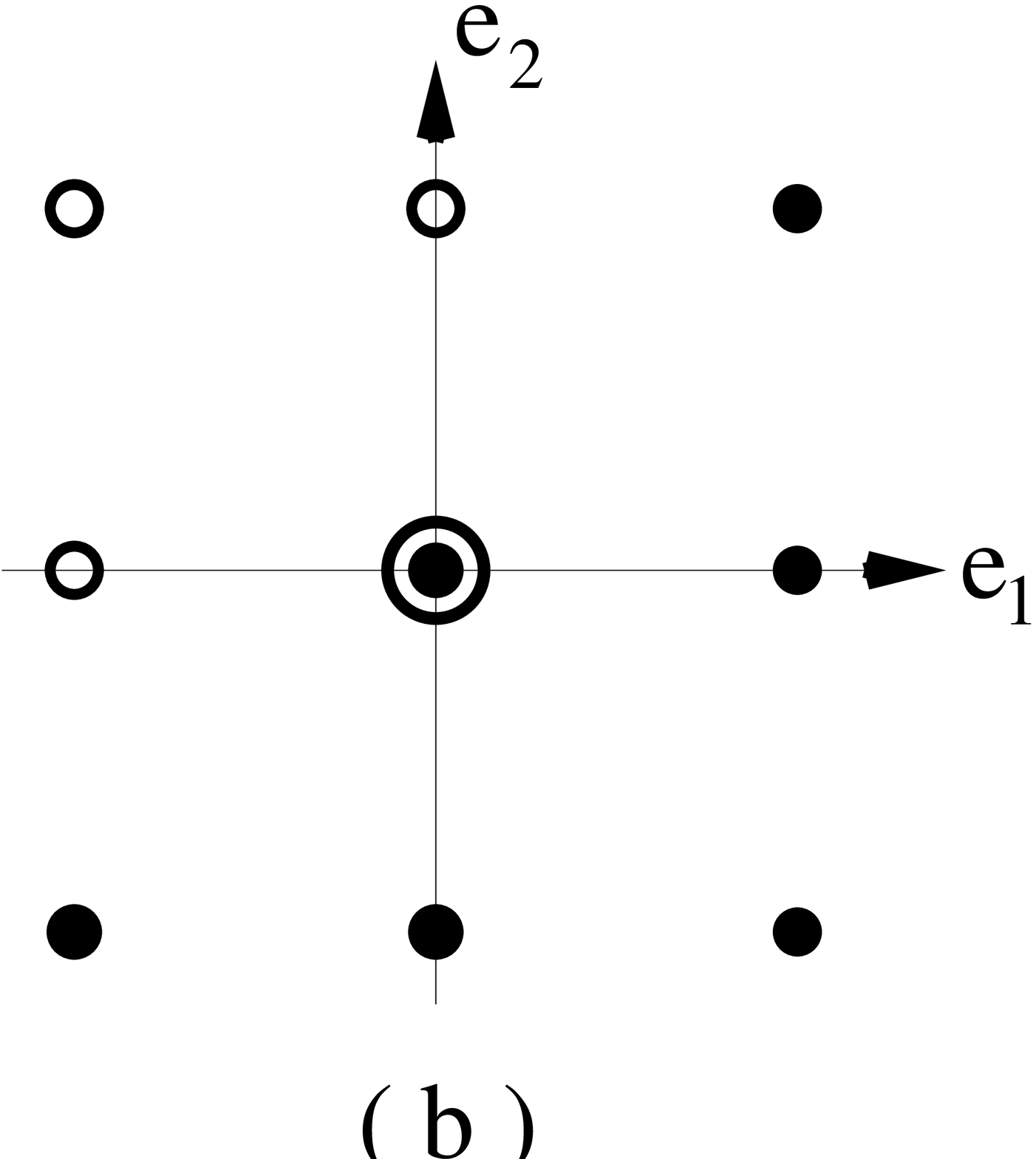} ~
\epsfxsize=1.5in\epsfbox{
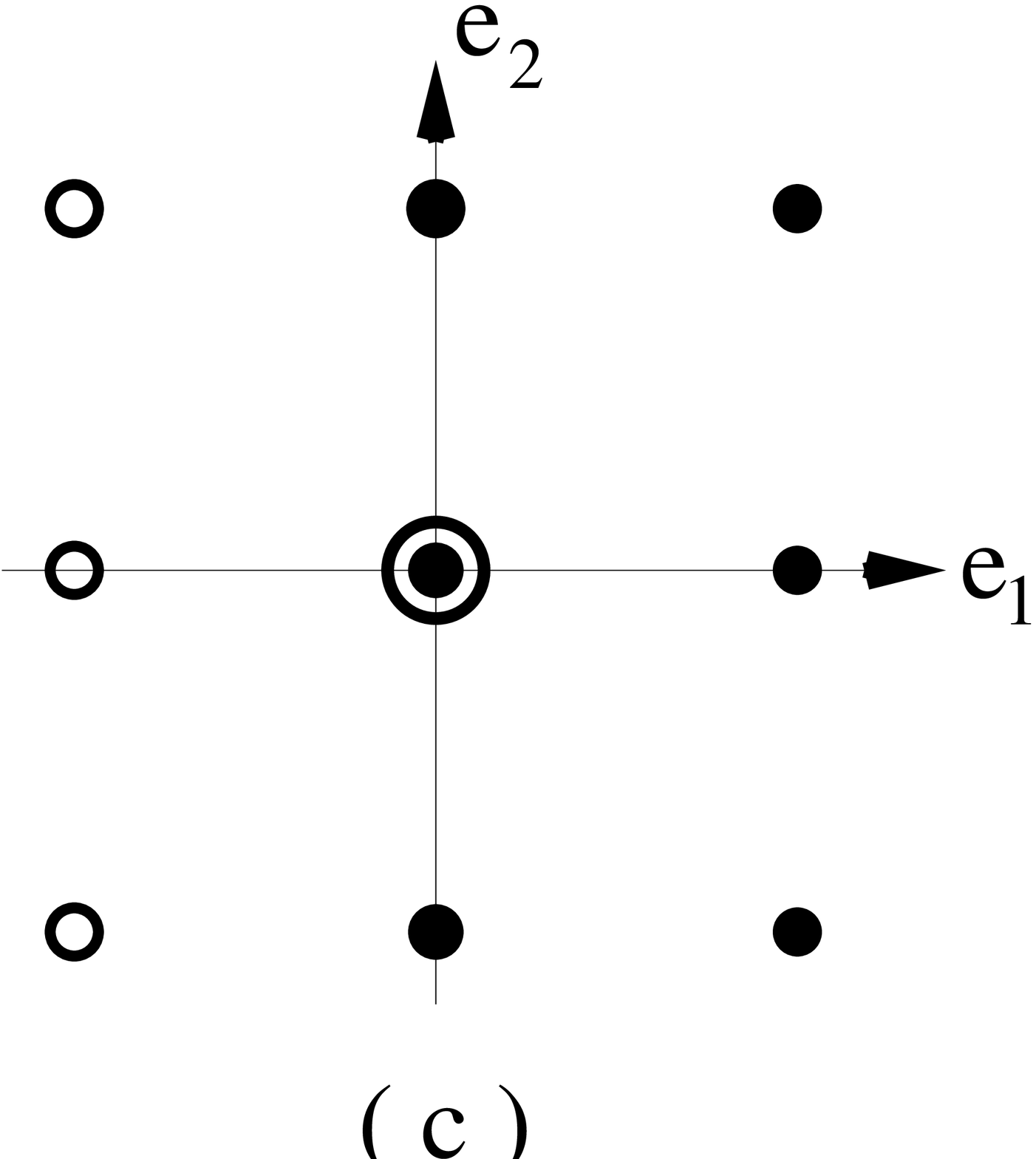}
}
\caption[Root space]{(a) Root diagram 
of the complexified conformal Lie algebra
$(\mathfrak{conf}_3)_{\mathbb{C}}$ and the labelling of its generators.
The double circle in the center denotes the Cartan subalgebra.   
The generators of the maximal parabolic subalgebras
(b) ${\mathfrak{sch}}$ and
(c) ${\mathfrak{alt}}$ are indicated by the filled dots. 
\label{Bild1}}
\end{figure}

Ageing phenomena occur widely in physics: glasses, granular systems or phase-ordering kinetics are just a few examples. While it is well-accepted that
they display some sort of dynamical scaling, the question has been raised whether their non-equilibrium dynamics might posses larger symmetries than merely scale-invariance. At first sight, the noisy terms in the Langevin
equations usually employed to model these systems might appear to exclude
any non-trivial answer, but it was understood recently that provided the
deterministic part of a Langevin equation is Galilei-invariant, then all
observables can be exactly expressed in terms of multipoint correlation functions calculable from the deterministic part only.\cite{Pico04} It is therefore of interest to study the dynamical symmetries of non-linear partial differential equations which extend dynamical scaling. In this context, the 
so-called Schr\"odinger algebra $\mathfrak{sch}$ has been shown to play an 
important r\^ole in phase-ordering kinetics. In what follows we shall restrict 
to one space dimension and we recall in figure~\ref{Bild1} through a root 
diagram the definition of $\sch$ as a parabolic subalgebra of the conformal 
algebra $\conf_3$.\cite{Henk03} 
%%----------------------------------------------------------------------------%%
%%----------------------------------------------------------------------------%%

\section{A brief perspective on the algebra $\alt$}
There is a classification of semi-linear partial differential equations with a 
parabolic subalgebra of $\conf_3$ as a symmetry.\cite{MS} Here 
we shall study the abstract Lie algebra $\alt$, which is the other maximal
parabolic subalgebra of $\conf_3$ (see figure~\ref{Bild1}) 
and its representations; we shall also see that, like the algebra $\sch$, it can be embedded naturally
in an infinite-dimensional Lie algebra ${\mathcal W}$ which is an 
extension of the algebra $\mathrm{Vect}(S^1)$ of vector fields on the circle. Quite strikingly, we shall find on our way a
'no-go theorem' that proves the impossibility of a conventional extension of the embedding $\alt\subset\conf_3$.

\subsection{The abstract Lie algebra $\alt$}

Elementary computations make it clear that (see figure~\ref{Bild1} and $D=2X_0-N$)
\BEQ
\alt=\langle V_+,D,Y_{-\half}\rangle\ltimes \langle X_1,Y_{\half},M_0\rangle:=\g\ltimes\h
\EEQ
is a semi-direct product of $\g\simeq \slin(2,\R)$ by a three-dimensional commutative Lie algebra
$\h$; the vector space $\h$ is the irreducible spin-1 real representation of $\slin(2,\R)$, which
can be identified with $\slin(2,\R)$ itself with the adjoint action.
So one has the following

\noindent {\bf Proposition 2.1:}

\begin{enumerate}
\item
$ \alt\simeq \slin(2,\R)\otimes\R[\eps]/\eps^2,$
{\it where $\eps$ is a 'Grassmann' variable;}
\item
$\alt\simeq\p_{3}$
{\it where $\p_3\simeq\so(2,1)\ltimes\R^3$ is the relativistic Poincar\'e algebra in (2+1)-dimensions.}
\end{enumerate}

\noindent {\bf Proof :}
The linear map $\Phi:\alt\to\slin(2,\R)\otimes\R[\eps]/\eps^2$ defined by
$$
\Phi(V_+)=L_1,\ \Phi(D)=L_0,\ \Phi(Y_{-\half})=L_{-1} 
$$
$$ 
\Phi(X_1)=\half L_1^{\eps}, \ \Phi(Y_{\half})=L_0^{\eps},\ \Phi(M_0)=L_{-1}^{\eps} 
$$
is easily checked to be a Lie isomorphism. \hfill \eop

In particular, the representations of $\alt\simeq\p_3$
are  well-known since Wigner studied them in the 30'es.

\subsection{Central extensions: an introduction}
Consider any Lie algebra $\g$ and an antisymmetric real two-form $\alpha$ on $\g$. 
Suppose that its Lie bracket $[\ , ]$ can be 'deformed' into a new
Lie bracket $\widetilde{[\ , ]}$ on $\tilde{\g}:=\g\times \R K$, where $[K,\g]=0$,
 by putting $\widetilde{[(X,0),(Y,0)]}=([X,Y],\alpha(X,Y))$. Then $\tilde{\g}$ is called a
{\it central extension} of $\g$. The Jacobi identity is equivalent with 
the nullity of the totally antisymmetric three-form 
$\D\al:\Lambda^3(\g)\to\R$ defined by 
$$
\D\alpha(X,Y,Z)=\al([X,Y],Z)+\al([Y,Z],X)+\al([Z,X],Y).
$$
Now we say that two central extensions $\g_1,\g_2$ of $\g$ defined by $\al_1,\al_2$ are equivalent if
$\g_2$ can be gotten from $\g_1$ by substituting $(X,c)\mapsto (X,c+\lambda(X))$ $(X\in\g)$ for
a certain 1-form $\lambda\in\g^*$, that is, by changing the non-intrinsic embedding of $\g$
into $\tilde{\g}_1$. In other words, $\al_1$ and $\al_2$ are equivalent if
$\al_2-\al_1=\D\lambda$, where $\D\lambda(X,Y)=\langle \lambda,[X,Y]\rangle.$ 
The operator $\D$ can be made into the differential of a complex (called Chevalley-Eilenberg complex),
and the preceding considerations make it clear that the classes of equivalence of
central extensions of $\g$ make up a vector space $H^2(\g)=Z^2(\g)/B^2(\g)$, where $Z^2$ is 
the space of {\it cocycles} $\alpha\in\Lambda^2(\g^*)$ verifying $\D\alpha=0 $,
and $B^2$ is the space of {\it coboundaries} $\D\lambda$, $\lambda\in\g^*$. 
We have the well-known

\noindent {\bf Proposition 2.2:}
{\it The Lie algebra $\alt$ has no non-trivial central extension:
$H^2(\alt)=0$.}

All this becomes very different when one embeds $\alt$ into an infinite-dimensional Lie algebra.

\subsection{Infinite-dimensional extension of $\alt$}

The Lie algebra Vect($S^1$) of vector fields on the circle has a long story
in mathematical physics. It was discovered by Virasoro\cite{Vira70} in the 
70'es that $\mathrm{Vect}(S^1)$ has a 
one-parameter family of central extensions which yield the 
so-called {\em Virasoro algebra} 
\BEQ
{\mathfrak{vir}} := \mathrm{Vect}(S^1)\oplus\R K=\langle (L_n)_{n\in\Z},K\rangle
\EEQ
with Lie brackets 
$$
[K,L_n]=0, \quad [L_n,L_m]=(n-m)L_{n+m}+\delta_{n+m,0}\  c\ n(n^2-1)K\quad (c\in\R)
$$
When $c=0$, one retrieves Vect($S^1$) by identifying the $(L_n)$ with the 
usual Fourier basis $(e^{\II n\theta}\D\theta)_{n\in\Z}$ of periodic vector 
fields on $[0,2\pi]$, or with $-z^{n+1}\frac{\D}{\D z}$ with $z:=e^{\II\theta}$. 
Note in particular that $\langle L_{-1},L_0,L_1\rangle$
is isomorphic to $\slin(2,\R)$, and that the Virasoro cocycle
restricted to $\slin(2,\R)$ is $0$, as should be (since $\slin(2,\R)$ has no non-trivial central extensions).

It is tempting to embed $\alt\simeq \slin(2,\R)\otimes\R[\eps]/\eps^2$ into the Lie algebra
\BEQ
{\mathcal W}:=\mathrm{Vect}(S^1)\otimes\R[\eps]/\eps^2=\langle L_n\rangle_{n\in\Z}\ltimes \langle L_n^{\eps}
\rangle_{n\in\Z},
\EEQ
with Lie brackets
$$ \label{gl:W}
[L_n,L_m]=(n-m)L_{n+m},\ [L_n,L_m^{\eps}]=(n-m)L_{n+m}^{\eps},\ [L_n^{\eps},L_m^{\eps}]=0.
$$
These brackets come out naturally putting $\mathcal W$ in the $2\times 2$-matrix form
\BEQ
L_n\mapsto \left( \begin{array}{cc} L_n & 0\\ 0& L_n\end{array}\right),
\ L_n^{\eps}\mapsto \left(\begin{array}{cc} 0& L_n\\ 0& 0 \end{array}\right)
\EEQ
leading to straightforward generalizations (there exists a deformation of
$\sv$ that can be represented as upper-triangular $3\times 3$ Virasoro 
matrices\cite{RogUnt05}).

In terms of the standard representations of Vect($S^1$) as modules of $\al$-densities
${\mathcal F}_{\al}=\{ u(z)(\D z)^{\al}\}$ with the action
\BEQ
f(z)\frac{\D}{\D z} \left(u(z)(dz)^{\al}\right)=(fu'+\al f' u)(z) (\D z)^{\al}, 
\EEQ 
we have

\noindent {\bf Proposition 2.3:}
${\mathcal W}\simeq {\mathrm{Vect}}(S^1)\ltimes {\mathcal F}_{-1}$.

There are two linearly independent central extensions of $\mathcal W$:

1. the natural extension to $\mathcal W$ of the Virasoro cocycle on Vect($S^1$), namely
$[\ ,\ ]=\widetilde{[\ ,\ ]}$ except for $\widetilde{[L_n,L_{-n}]}=n(n^2-1)K+2nL_0.$ In other
words, Vect($S^1$) is centrally extended, but its action on ${\mathcal F}_{-1}$ remains unchanged;

2. the cocycle $\omega$ which is zero on $\Lambda^2($Vect($S^1$)) and
$\Lambda^2({\mathcal F}_{-1})$, and defined by $
\omega(L_n,L_m^{\eps})=\delta_{n+m,0}\, n(n^2-1)K $
on Vect($S^1)\times {\mathcal F}_{-1}$. 

A natural related question is: can one deform the extension of Vect($S^1$) by 
the Vect($S^1$)-module ${\mathcal F}_{-1}$ ? The answer is: no, thanks to the triviality\cite{Fuks}\cite[chapter 4]{GuiRog05} of the cohomology space $H^2({\rm Vect}(S^1),{\mathcal F}_{-1})$.
Hence, any Lie algebra structure $\widetilde{[\ ,\ ]}$ on the vector space
Vect($S^1$)$\oplus{\mathcal F}_{-1}$ such that 
$$
\widetilde{[(X,\phi),(Y,\psi)]}=\left([X,Y]_{{\rm Vect}(S^1)}, 
{\rm ad}_{{\rm Vect}(S^1)}X.\psi-{\rm ad}_{{\rm Vect}(S^1)}Y.\phi
+B(X,Y)\right)
$$
is isomorphic to the Lie structure of $\mathcal W$ (where $B$ is an antisymmetric two-form on Vect($S^1$)). 
So one may say that $\mathcal W$ and its central extensions are natural objects to look at.
\subsection{Some results on representations of ${\mathcal W}$}
We now state two results which may deserve deeper thoughts and will be developed in the future.\\
\noindent {\bf Proposition 2.4 ({\rm 'no-go theorem'}):}
{\it There is no way to extend the usual representation of $\alt$ as conformal vector fields into
an embedding of $\mathcal W$ into the Lie algebra of vector fields on $\R^3$.}\\
\noindent {\bf Proposition 2.5:} {\it The infinite-dimensional extension ${\mathcal W}$ 
of the algebra $\mathfrak{alt}_1$ is a contraction of a pair 
of commuting Virasoro algebras 
$\mathfrak{vir}\oplus\overline{\mathfrak{vir}}\rar{\mathcal W}$. In particular, we
have the explicit differential operator representation}
\BEA
L_n &=& -t^{n+1}\partial_t +(n+1)t^n r\partial_r -(n+1)x t^n 
-n(n+1)\gamma t^{n-1} r \nonumber \\
L_n^{\eps} &=& -t^{n+1}\partial_r -(n+1)\gamma t^n 
\label{gl:formeexpl}
\EEA
{\it where $x$ and $\gamma$ are parameters and $n\in\mathbb{Z}$.}

\section{Appell systems}
\begin{definition}
{\em Appell polynomials} $\{h_n(x);n\in \mathbb{N}  \}$ on $\mathbb{R}$ are usually characterized by the two conditions
\begin{itemize}
\item $h_n(x)$ are polynomials of degree $n$,
\item $D h_n(x)=nh_{n-1}(x) $, where $D$ is the usual derivation opertor.
\end{itemize}
\end{definition}
Interesting examples are furnished by the shifted moment sequences
$$h_n(x)=\int_{-\infty}^{\infty}(x+y)^{n}p(\D y)$$
where $p$ is a probability measure on $\R$ with all moments finite.\\
This definition generalizes to higher dimensions.
On non-commutative algebraic structures, the shifting corresponds to left or right multiplication and in general, $\{h_n\}$ is not a family of polynomials. We shall call it {\em Appell systems}.\cite{FS}\\
Appell systems of the Schr\"odinger algebra $\sch$ have been 
investigated\cite{FKS} but the algebra $\alt$ requires a specific study.
$\alt$ has the following Cartan decomposition:
\begin{equation}\label{cartan}
\mathfrak{\alt}=\mathfrak{P}\oplus \mathfrak{K}\oplus \mathfrak{L} 
= \{ Y_1,X_1\} \oplus \{ Y_0,X_0\} \oplus \{ Y_{-1},X_{-1}\}
\end{equation}
and there is a one to one correspondence between the subalgebras $\mathfrak{P}$ and $\mathfrak{L}$.
Write $X\in \alt$ in the form 
$X=\sum_{i=1}^{6}{\alpha}_{i}b_i$, where $\{b_i, i=1,\ldots,6 \}$ is a basis of 
$\alt$. The $\alpha_{i}$ are called coordinates of the first kind.\\
Here we use the basis $b_1=Y_1$, $b_2=X_1$, $b_3=Y_0$, $b_4=X_0$, $b_5=Y_{-1}$, $b_6=Y_1$. 
Let ALT be the simply connected Lie group corresponding to  $\alt$. Group elements in a neighbourhood of the identity can be expressed as
$$e^{X}=e^{A_1b_1}\ldots e^{A_6b_6}$$
The $A_i$ are called {\em coordinates of the second kind}.\\
%The general theory\cite{FS} tells us that the $A_i$ are solutions of 
%a system of differential equations and that they contain the complete 
%information about the Lie algebra $\alt_1$. Computations (which we omit% here %due to the lack of space) lead to: $A_1=$,$A_2=$, $A_3=$, $A_4=$, 
%$A_5=$,$A_6=$, 
Referring to the decomposition (\ref{cartan}), we specialize variables, writing $V_1,V_2,B_1,B_2$ for $A_1,A_2,A_5,A_6$ respectively. 
Basic for our approach is to establish the partial group law:
$e^{B_1Y_{-1}+B_2X_{-1}}e^{V_1Y_1+V_2X_1}=?$.
We get
\begin{equation}
B_1Y_{-1}+B_2X_{-1}=\left(\begin{array}{cccc} 0& 0& 0& 0\\ -B_2& 0& -B_1& 0\\ 0& 0& 0& 0\\ 0& 0& -B_2& 0\end{array}\right)
\;\; , \;\; 
V_1Y_1+V_2X_1=\left(\begin{array}{cccc} 0& V_2& 0& V_1\\ 0& 0& 0& 0\\ 0& 0& 0& V_2\\ 0& 0& 0& 0\end{array}\right)
\end{equation}
and finally:
\begin{equation}
e^{B_1Y_{-1}+B_2X_{-1}}e^{V_1Y_1+V_2X_1}=\left(\begin{array}{cccc} 1& V_2& 0& V_1\\ -B_2& 1-B_2V_2& -B_1& -B_2V_1-B_1V_2\\ 
0& 0& 1& V_2\\ 0& 0& -B_2& 1-B_2V_2\end{array}\right)
\end{equation}
\begin{proposition}
In coordinates of the second kind, we have the Leibnitz formula,
\begin{eqnarray}\label{eq:pglaw}
& & g(0,0,0,0,B_1,B_2)g(V_1,V_2,0,0,0,0)= g(A_1,A_2,A_3,A_4,A_5,A_6)=\nonumber \\
& & g({B_1V_2^2+V_1\over (1-B_2V_2)}, {V_2\over (1-B_2V_2)}, -2{B_1V_2+B_2V_1\over (1-B_2V_2)},
2\ln (1-B_2V_2),\\ 
& & {B_1-2B_1B_2V_2-B_2^2V_1\over (1-B_2V_2)^2},\nonumber  
 {B_2\over (1-B_2V_2)}) 
\end{eqnarray}
\end{proposition}
Now we are ready to construct the representation space and basis-the canonical Appell system. To start, 
define a vacuum state $\Omega$. The elements $Y_1,X_1$ of $\mathfrak{P}$
can be used to form basis elements
\begin{equation}
|jk\rangle =Y_1^jX_1^k\Omega ,j,k\geqslant 0
\end{equation}
of a Fock space $\mathfrak{F}=\mbox{\rm span}\{|jk\rangle \}$ on which $Y_1,X_1$ act as raising
operators, $Y_{-1},X_{-1}$ as lowering operator and $Y_0,X_0$ as multiplication
with the constants $\gamma,x $ (up to the sign) correspondingly. That is,
\begin{eqnarray}
&&Y_1\Omega =|10\rangle , X_1\Omega=|01\rangle \nonumber \\
&&Y_{-1}\Omega=0, X_{-1}\Omega=0 \label{eq:fock} \\
&&Y_0\Omega=-\gamma |00\rangle , X_0\Omega = -x|00\rangle \nonumber
\end{eqnarray}
The goal is to find an abelian subalgebra spanned by some selfadjoint
operators acting on representation space, just constructed. Such a two-dimensional subalgebra can be obtained by an appropriate ``turn'' of the plane $\mathfrak{P}$ in the Lie algebra, namely via the adjoint action of the group element formed by exponentiating $X_{-1}$. The resulting plane, $\mathfrak{P}_{\beta}$ say, is abelian and is spanned by
\begin{eqnarray}
\bar{Y_1}=e^{\beta X_{-1}}Y_1e^{-\beta X_{-1}}=Y_1-2\beta Y_0 +\beta^2Y_{-1} \nonumber \\
\bar{X_1}=e^{\beta X_{-1}}X_1e^{-\beta X_{-1}}=X_1-2\beta X_0 +\beta^2X_{-1} \label{eq:abepl} 
\end{eqnarray}
Next we determine our canonical Appell systems. We apply the Leibniz formula (\ref{eq:pglaw}) with 
$B_1=0, B_2=\beta ,V_1=z_1, V_2=z_2$ and (\ref{eq:fock}). This yields   
\begin{eqnarray}\label{eq:genfu}
e^{z_1\bar{Y_1}}e^{z_2\bar{X_1}}\Omega &=& e^{\beta X_{-1}}e^{z_1Y_1}e^{z_2X_1}e^{-\beta X_{-1}}\Omega=
e^{\beta X_{-1}}e^{z_1Y_1}e^{z_2X_1}\Omega= \nonumber \\
&=& e^{z_1Y_1\over (1-\beta z_2)^2}e^{z_2X_1\over (1-\beta z_2)}e^{2\gamma\beta z_1\over (1-\beta z_2)}(1-\beta z_2)^{-2x}\Omega
\end{eqnarray}
To get the generating function for the basis $|jk\rangle$ set in equation (\ref{eq:genfu})
\begin{equation}\label{eq:chang}
v_1={z_1\over (1-\beta z_2)^2},  v_2={z_2\over (1-\beta z_2)}
\end{equation}
Substituting throughout, we have
\begin{proposition}
The generating function for the canonical Appell system, $|jk\rangle =Y_1^jX_1^k\Omega$ is
\begin{eqnarray}\label{eq:appsy}
e^{v_1Y_1+v_2X_1}\Omega 
&=& exp(y_1{v_1\over (1+\beta v_2)^2})exp(y_2{v_2\over (1+\beta v_2)})
exp(-{2\gamma \beta v_1\over (1-\beta v_2)})\nonumber \\
& & (1+\beta v_2)^{-2x}\Omega
\end{eqnarray}
where we identify $\bar{Y_1}\Omega =y_1 \cdot 1$ and $\bar{X_1}\Omega =y_2 \cdot 1$
in the realization as function of $y_1,y_2$.
\end{proposition}
With $v_1=0$, we recognize the generating function for the Laguerre polynomials, while $v_2=0$ reduces to the generating function for Hermite polynomials.
%\section{Leibniz fuction and coherent states}

%\section{Concluding remarks}


\begin{thebibliography}{99}
\bibitem{FKS} P. Feinsilver, Y. Kocik and R. Schott, 
{\it Representations of the Schr\"odinger Algebra and Appell Systems}, 
Progress of Physics, {\bf 52} (2004) 343-359.

\bibitem{FS} P. Feinsilver and R. Schott, {\it Algebraic Structures and Operator Calculus, Vol.3: Representations of Lie Groups.} Kluwer Academic Publishers, Dordrecht, 1993.

\bibitem{Fuks} D. B. Fuks, {\it Cohohomology of infinite-dimensional Lie algebras}, Contemporary Soviet Mathematics, Consultants Bureau, New York, 1986.

\bibitem{MSSU} M. Henkel, R. Schott, S. Stoimenov and J. Unterberger, 
{\it On the dynamical symmetric algebra of ageing: Lie structure, representations and Appell systems.} Pr\'epublication Institut Elie Cartan, 2005.

\bibitem{GuiRog05} L. Guieu, C. Roger, preprint, available on \\ {\tt http://www.math.univ-montp2.fr/\~{}guieu/The\_Virasoro\_Project/Phase1/}.

\bibitem{Henk03} M. Henkel, J. Unterberger, {\it Schr\"odinger-invariance and space-time symmetries}, Nuclear Physics {\bf B660} (2003) 407-435.

\bibitem{Pico04} A. Picone and M. Henkel, {\it Local scale-invariance and ageing in noisy systems}, Nuclear Physics {\bf B688} (2004) 217-265. 

\bibitem{RogUnt05} C. Roger and J. Unterberger, in progress (2005). 

\bibitem{MS} S. Stoimenov and M. Henkel, 
{\it Dynamical symmetries of semi-linear Schr\"odinger and diffusion equations}, 
Nuclear Physics {\bf B723} (2005) 205-233.

\bibitem{Vira70} M. Virasoro, Phys. Rev. {\bf D1}, 2933-2936 (1970).


\end{thebibliography}
\end{document}